\documentclass{amsart}
\usepackage[dvips]{graphicx}

\newtheorem{thm}{Theorem} 
\newtheorem{lem}[thm]{Lemma}
 
\newtheorem{conj}[thm]{Conjecture}

\theoremstyle{definition}
\newtheorem{defn}{Definition}
\newtheorem{rem}{Remark}
\newtheorem{exmp}{Example}

\begin{document}

\title[Artin Groups of Large Type]
	{Three-generator Artin Groups of Large Type are Biautomatic}

\author[T.~Brady]{Thomas Brady}

	\address{Dept. of Mathematics\\
		Dublin City University\\
		Glasnevin, Dublin 9\\
		Ireland}
	\email{tom.brady@dcu.ie}

\author[J.~M{\tiny c}Cammond]{Jonathan P. M{\tiny c}Cammond}

	\address{Dept. of Mathematics\\ 
		Texas A\&M University\\ 
		College Station, TX 77843\\
		USA}
	\email{jon.mccammond@math.tamu.edu} 



\date{\today}

\begin{abstract}

In this article we construct a piecewise Euclidean, non-positively
curved $2$-complex for the $3$-generator Artin groups of large type.
As a consequence we show that these groups are biautomatic.  A slight
modification of the proof shows that many other Artin groups are also
biautomatic.  The general question (whether all Artin groups are
biautomatic) remains open.

\end{abstract}

\maketitle

\section{Introduction}

In this article we construct a piecewise Euclidean, non-positively
curved $2$-complex for the $3$-generator Artin groups of large
type.  Although standard presentations for the Artin groups typically
use words which alternate between pairs of the generators, there are
other presentations with `nicer' properties.  We shall be employing
one of these.

If we let $(a,b)_k$ denote the word of the form $ababa\ldots$ with
exactly $k$ letters, then the standard presentation of an Artin group
with three generators is given as 
\begin{equation}\label{stnd-pres}
G_{m,n,p} = \langle \, a,b,c \mid (a,b)_m = (b,a)_m,\ (b,c)_n = (c,b)_n,\
(c,a)_p = (a,c)_p \, \rangle
\end{equation}

\noindent
where $m$, $n$, and $p$ are positive integers greater than $1$.  After
introducing an alternative presentation, we will show that the new
presentation satisfies the small cancellation conditions $C(3)-T(6)$.
As a consequence, a metric of non-positive curvature can be defined on
the standard $2$-complex of this presentation and by \cite{gersten} the
group it defines is biautomatic.  In the final section, we present
a generalization of this result.

\section{The Presentation}

The non-positively curved $2$-complex alluded to above is constructed
by altering the standard presentations for Artin groups.  Consider the
following families of group presentations.
\begin{equation}
G_m = \langle \, a_1, a_2 \mid (a_1 ,a_2)_m = (a_2, a_1)_m \,\rangle
\end{equation}
\begin{equation}
H_{2k} = \langle \, x, a_1 \mid x^ka_1 = a_1x^k \rangle
\end{equation}
\begin{equation}
H_{2k+1} = \langle \, x, a_1 \mid x^{k+1} = a_1x^ka_1 \rangle
\end{equation}
\begin{equation}
I_m = \langle \, x, a_1, a_2, \ldots a_m \mid x = a_1a_2, x=a_2a_3, \ldots, x=a_ma_1 \rangle
\end{equation}

\noindent
Note that $G_m$ is the standard presentation for Artin groups with two
generators.  We will now show that the other presentations define the
same group.

\begin{lem}\label{equivalent}

For all $m > 1$, $G_m \cong H_m \cong I_m$.

\end{lem}

\begin{figure}
  \begin{center}~
  \includegraphics{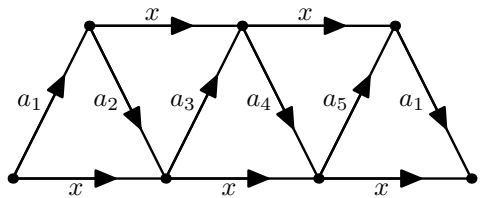}
  \end{center}
  \caption{Relations for $H_5$ and $I_5$
  \label{relations}}
\end{figure}

\begin{proof}

Start with $G_m$. If we add a new generator $x$ to $G_m$ with
definition $x = a_1a_2$, and then use this relation to eliminate $a_2$
we will end up with the presentation $H_m$.  The precise form of the
presentation will depend on whether $m$ is even or odd.  Since this
process can also be reversed, $G_m \cong H_m$.

The relationship between $H_m$ and $I_m$ is even more
straight-forward.  The relations of $I_m$ are derived by partitioning
the unique relation of $H_m$, while the relation in $H_m$ is formed by
arranging the relations of $I_m$ to eliminate the generators $a_2$,
$a_3$, \ldots, $a_m$.  This correspondence is illustrated in
Figure~\ref{relations}.
\end{proof}

\begin{lem}\label{presentation}

The Artin group $G_{m,n,p}$ can be presented as follows:
\[
I_{m,n,p} = 
\left<
\begin{array}{c}
x, a, d_3, \dots d_m\\
y, b, e_3, \dots , e_n\\ 
z, c, f_3, \dots , f_p 
\end{array}
\right.
\left|
\begin{array}{c}
x = ab,\, x = bd_3,\, x = d_3d_4, \dots, x = d_ma\\ 
y = bc,\, y = ce_3,\, y = e_3e_4, \dots, y = e_nb \\
z = ca,\, z = af_3,\, z = f_3f_4, \dots, z = f_pc
\end{array}
\right>
\]
\end{lem}

\begin{proof}

If the Tietze transformations described in the proof of the previous
lemma are applied to each of the sets of relations in this
presentation, the result will be the standard presentation
(Equation~\ref{stnd-pres}).  Notice that $a$ and $b$ play the role of
$d_1$ and $d_2$ in the first set of equations, $b$ and $c$ play the
role of $e_1$ and $e_2$ in the second set, and $c$ and $a$ play the
role of $f_1$ and $f_2$ in the third set.
\end{proof}

\section{The $2$-Complex}

Let $K_{m,n,p}$ be the standard $2$-complex associated with the
presentation $I_{m,n,p}$ given in Lemma~\ref{presentation}.  The
complex $K_{m,n,p}$ has a unique $0$-cell, $m+n+p$ $1$-cells (one for each of
the generators), and $m+n+p$ triangular $2$-cells (one for each of the
relations).  Of particular interest is the link the unique $0$-cell.
Recall that the link of a $0$-cell in a $2$-complex is the graph
obtained by intersecting a small sphere centered at the $0$-cell with
the surrounding complex.  In this case, each $1$-cell will contribute
two vertices to the link, and each corner of a $2$-cell will
contribute an edge.  The link of $K_{m,n,p}$ will be refered to as
$L_{m,n,p}$.

\begin{figure}
  \begin{center}
  \includegraphics{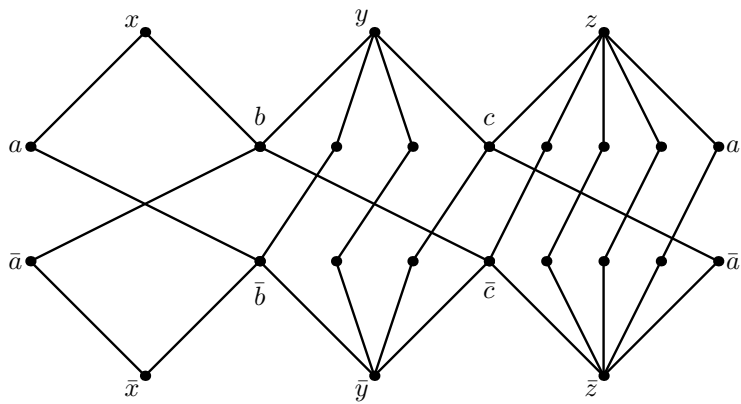}
  \end{center}
  \caption{The link of the vertex in $G_{2,4,5}$
  \label{example}}
\end{figure}

\begin{exmp}

The link $L_{2,4,5}$ is shown in Figure~\ref{example}, with the
exception that the two vertices on the extreme right need to be
identified with the two vertices on the extreme left.  We have adopted
the convention that letters without bars denote the vertex contributed
by the head of the oriented $1$-cell and letters with bars denote the
vertex contributed by the tail.

\end{exmp}

\begin{rem}

Notice that the link is essentially composed of three links of
$2$-generator Artin groups joined together along pairs of vertices.
In the example shown, we see the link of $I_2$, $I_4$, and $I_5$ laid
side by side.  Each of these links is in turn made up of packets of
paths of length $3$ which are twisted together.  The subgraphs
corresponding to the $2$-generator Artin groups will be called the
local pieces of the link.  Notice also that these local pieces overlap
only at vertices which arise from the standard generators $a$, $b$,
and $c$.  
Another key property of the link is that the vertices separate
themselves into $4$ distinct levels, with $\bar{x}$, $\bar{y}$, and
$\bar{z}$ at the bottom (level $1$), and $x$, $y$, and $z$ at the top
(level $4$).  The edges connecting levels $2$ and $3$ will be called
middle edges.

\end{rem}

In Figure~\ref{example} there are two short loops of length $4$.
Namely the paths connecting $\bar{a}-b-y-c-\bar{a}$ and
$\bar{a}-b-\bar{c}-\bar{z}-\bar{a}$.  We will now show that if there
are no commutation relations, then no short embedded loops can exist.

\begin{lem}\label{short loops}

If $m$, $n$, and $p$ are all at least $3$ then every embedded loop in
$L_{m,n,p}$ contains at least $6$ edges.

\end{lem}

\begin{figure}
  \begin{center}~
  \includegraphics{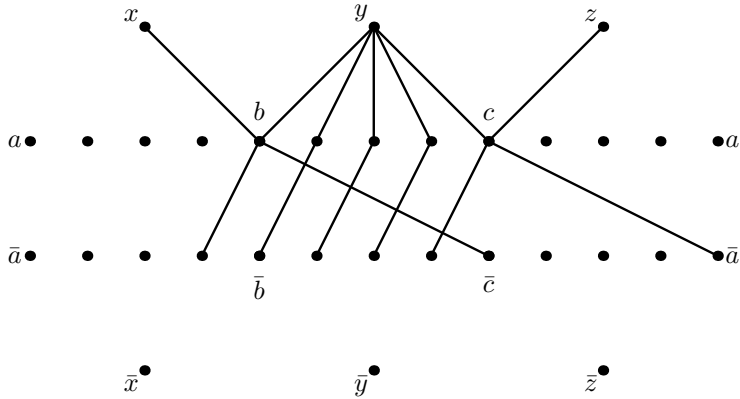}
  \end{center}
  \caption{The radius $2$ neighborhood of $y$ in the link $L_{5,5,5}$
  \label{neighborhood}}
\end{figure}

\begin{proof}

(Sketch) Since every edge in the link joins a pair of vertices on
adjacient levels, we see that the link is bipartite, and thus every
embedded loop will have an even length.  Furthermore, it should be
clear from the example that embedded loops of length $2$ cannot occur,
so that it only remains to examine possible loops of length $4$.

If there was a loop of length $4$ which contained a top vertex or a
bottom vertex, then the entire loop would be contained in the radius
two neighbourhood of this vertex.  An examination however shows this
neighborhood is a tree.  As an example, the radius two neighbourhood
of $y$ in $L_{5,5,5}$ is shown in Figure~\ref{neighborhood}.
On the otherhand, if there was such a loop which did not pass through
a top or bottom vertex then it would be completely composed of middle
edges.  The subcomplex of these edges however consists of $m+n+p-9$
disjoint edges and three chains of length $3$.  This completes the
proof.
\end{proof}

In the language of small cancellation theory, this lemma shows that
the presentation $I_{m,n,p}$ satisfies the condition $T(6)$.

\section{The Main Result}

The complex $K_{m,n,p}$ can be given a piecewise Euclidean metric by
assigning length of $1$ to each $1$-cell and making each $2$-cell an
equilateral triangle.  This metric also induces a metric on the link
of the $0$-cell by assigning the angle at a corner of a $2$-cell as
the length of the corresponding edge in the link.  In this instance,
the triangles are all equilateral.  Thus all of the edges in the link
are assigned a length of $\frac{\pi}{3}$.

\begin{thm}

If $m$, $n$, and $p$ are all at least $3$, then the metric on
$K_{m,n,p}$ is non-positive curved and the group $G_{m,n,p}$ is
biautomatic.  In particular, all $3$-generator Artin groups of large
type are biautomatic.

\end{thm}

\begin{proof}

By Lemma~\ref{short loops}, every embedded loop in the link of the
$0$-cell contains at least six edges, and thus it measures at least
$2\pi$.  By \cite{ballmann} this proves that the metric on $K_{m,n,p}$
is non-positively curved.
That the group is biautomatic follows from a result of Gersten and
Short (\cite{gersten}).  In their terminology $K_{m,n,p}$ is an $A_2$
complex.
\end{proof}

\section{Generalizations}

In order to extend this result Artin groups with more than three
generators, additional terminology needs to be introduced.

\begin{defn}[$G_\Gamma$]

An arbitrary Artin group can be defined as follows.  Let $\Gamma$ be a
undirected graph with no loops or multiple edges and with a positive
integer greater than $1$ assigned to each edge.  The standard
presentation for the Artin group associated to $\Gamma$ is defined as
follows: it has one generator for each of the vertices of $\Gamma$ and
for each edge there is a relation.  In particular, if there is an edge
connected a vertex $a$ to a vertex $b$ which is labeled by $m$ then we
add the relation $(a,b)_m = (b,a)_m$.  We will call this presentation
$G_\Gamma$.

\end{defn}

\begin{defn}[$I_\Gamma$, $K_\Gamma$, and $L_\Gamma$]

In order to convert this presentation into a presentation like that
used in Lemma~\ref{presentation}, we need to include additional
information.  There is an ambiguity which stems from the fact that the
presentation $G_m$ is symmetric with respect to the generators $a_1$
and $a_2$, but the equivalent presentation $I_m$ is not.  In
particular, notice that in $I_m$, $a_1a_2$ occurs as a subword of a
relation, but that $a_2a_1$ does not.  Thus to precisely define a
presentation of the type we have been considering, we need to require
that the defining graph $\Gamma$ is a directed graph.  The direction
will indicate which subword of length $2$ to be preserved in the
rewriting process.  The presentation which results will be called
$I_\Gamma$, the standard $2$-complex for this presentation will be
denoted $K_\Gamma$, and the link of its unique $0$-cell will be
denoted $L_\Gamma$.
\end{defn}

\begin{figure}
  \begin{center}~ 
  \includegraphics{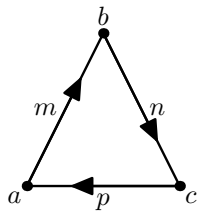} 
  \end{center}
  \caption{Defining graph $\Gamma$ for the presentation $I_{m,n,p}$
  \label{mnp}}
\end{figure}

\begin{exmp}

In Figure~\ref{mnp} we have displayed the defining graph $\Gamma$
which leads to the presentation given in Lemma~\ref{presentation}.
Although it was not highlighted at the time, the orientations of the
edges of $\Gamma$ were in fact crucial to the argument in the proof of
Lemma~\ref{short loops}.  As we will see in Lemma~\ref{graphs}, if one
of the orientations had been reversed then embedded loops of length
$4$ would indeed have been present in the link.

\end{exmp}

\begin{rem}

As in the $3$-generator case the link $L_\Gamma$ is a bipartite graph
with $4$ distinct levels of vertices, and since it is bipartite, all
of its embedded loops will have an even length.
The link itself is comprised of various local pieces which look like
the link of $I_m$ for some $m$.  Moreover, any two of these local
pieces are joined to each other only at a pair of vertices which
correspond to a generator of $G_\Gamma$ that arises from a vertex of
$\Gamma$.  For example, in Figure~\ref{example} the local pieces which
look like the links of $I_4$ and $I_5$ are joined together at the
vertices $c$ and $\bar{c}$, and these vertices correspond the
generator $c$ in $G_\Gamma$ that corresponds to the vertex labeled $c$
in Figure~\ref{mnp}.  The vertices in $L_\Gamma$ which correspond in
this way to the vertices of $\Gamma$ will be called special vertices
and the edges whose endpoints are both special will be called special
edges.

\end{rem}

\begin{rem}

There is another connection between $\Gamma$ and $L_\Gamma$ which will
be useful in the proof given below.  If the vertices of the link are
organized so that the vertices fall into four distinct levels, the
pairs of vertices are aligned to that the one with the bar is directly
below the other, and the local pieces of the link lie in a vertical
$2$-dimensional surface, then the graph $\Gamma$ can be thought of as
the ``top-view'' of the link.  In this view the pairs of special
vertices where the local pieces overlap are seen as the vertices of
$\Gamma$, the local pieces which are in one-to-one correspondence with
the edges of $\Gamma$ are seen as being projected onto the appropriate
edges, and the special edges in the link which allow one to get from a
special vertex in level $2$ to a special vertex in level $3$ in a
single step, are seen as the ones which give the orientation to the
edge corresponding to its local piece.

\end{rem}

The lemma about short embedded loops can now be extended to the
following, more general, result.

\begin{figure}
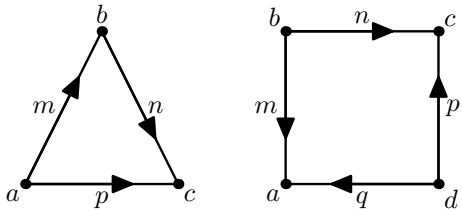

  \begin{center}~
  \includegraphics{badgraph-a.eps}
  \hspace{.2in}
  \includegraphics{badgraph-b.eps}
  \end{center}
  \caption{Subgraphs which lead to short embedded loops
  \label{bad graphs}}
\end{figure}

\begin{lem}\label{graphs}

The link $L_\Gamma$ contains an embedded loop of length $4$ if and
only if the defining graph $\Gamma$ contains an oriented subgraph
isomorphic to one of the graphs in Figure~\ref{bad graphs}.  Moreover,
these loops of length $4$ are the only possible embedded loops of
length less than $6$.

\end{lem}

\begin{proof}

We will begin by showing that all short loops must involve at least
$3$ local pieces.  Clearly the local pieces themselves have no
embedded loops of length less than $6$, so consider a subgraph of the
link consisting of only $2$ (connected) local pieces and let $a$ and
$\bar{a}$ be the pair of vertices connecting them.  Since the shortest
path connecting $a$ to $\bar{a}$ has length $3$ and the shortest
embedded path connecting either vertex to itself has length $6$, all
of the loops in this subgraph will have length at least $6$.

Next, we will establish the forward direction.  Assume that $L_\Gamma$
contains a short embedded loop of length $4$.  If this loop does not
contain a top or a bottom vertex, then it is comprised entirely of
middle edges.  Since the special vertices are the only vertices in
levels $2$ or $3$ which are connected to more than one middle edge,
all four vertices in the loop must be special and as a result all four
of the edges are special as well.  The vertices which are adjacient in
the loop must correspond to distinct vertices in $\Gamma$ since they
connected by a special edge, and the vertices which are not adjacient
in the loop must also correspond to distinct vertices in $\Gamma$
since the loop is embedded.  This shows that the subgraph of $\Gamma$
determined by these four special vertices contains the graph on the
righthand of Figure~\ref{bad graphs} as a subgraph.

If this short loop contained more than one top or bottom vertex then
it would be contained in the union of two local pieces since all of
the edges connected to a top or bottom vertex belong to the same local
piece.  We have already shown that the union of two local pieces do
not contain short loops, so there are no short loops of this type.

Finally, suppose this short loop contained exactly one top vertex.
Since the edges connected to it belong to the same local piece and
since we know that at least three local pieces are involved, the other
two edges must lie in distinct local pieces.  Thus all three of the
other vertices are special, and the two edges not connected to the top
vertex are special.  Since three distinct local pieces are involved it
is clear that the three special vertices correspond to distinct
vertices in $\Gamma$.  At this point we know that the subgraph of
$\Gamma$ determined by these these three local pieces is a triangle,
and we know the orientations of two of the sides.  The orientation of
the third side (corresponding to the local piece containing the top
vertex) is unknown, but with either orientation, the resulting graph
is isomorphic to the graph on the lefthandside of Figure~\ref{bad
graphs}.  Since the case of exactly one bottom vertex is similar, the
proof in the forward direction is complete.  The proof in the other
direction is essentially immediate.  If either of these graphs appears
as a subgraph of $\Gamma$, then it is easy to exhibit a path of length
$4$ in the link.
\end{proof}

The embedded loops of length $4$ which result from the graph on the
left will be called loops of type $A$ while those which result from
the one on the right will be called loops of type $B$.

\begin{thm}\label{triangle-free}

If $\Gamma$ does not contain any triangles, then $K_\Gamma$ has a
metric of non-positive curvature and the Artin group $G_\Gamma$ is
biautomatic.

\end{thm}

\begin{proof}

Since by hypothesis loops of type $A$ cannot occur, we can conclude
that all loops of length $4$ contain $4$ middle edges (since this is
true for loops of type $B$).  Next, we claim that all loops of length
$6$ must contain at least $2$ middle edges.  The structure of the link
guarantees that a loop will contain an even number of middle edges.
If it did not contain any middle edges, then it would be contained in
either the top or the bottom.  Since top and bottom vertices are only
connected to edges in the same local piece of the link, the loop
itself would be contained in three local pieces.  This would implies
that a triangle exists in $\Gamma$, contradiction.

To complete the proof we assign a metric to $K_\Gamma$ as follows.
The $1$-cells corresponding to the top and bottom vertices are
assigned a length of $\sqrt{2}$ and all others are assigned a length
of $1$.  The $2$-cells are Euclidean triangles with angles of
$\frac{\pi}{2}$, $\frac{\pi}{4}$, and $\frac{\pi}{4}$.  The conditions
shown in the previous paragraph, demonstrate that all embedded loops
in the link have a length of at least $2\pi$, and by \cite{ballmann}
the metric on $K_\Gamma$ is non-positively curved.  Gersten and Short
(\cite{gersten2}) again allow us to conclude that the group $G_\Gamma$
is biautomatic.  In their terminology $K$ is a complex of type $B_2$.
\end{proof}

This particular result was in some sense already known since Pride
showed that the standard presentation of these `triangle-free' Artin
groups satisfies the $C(4)-T(4)$ conditions (\cite{pride}), and
Gersten-Short showed that $C(4)-T(4)$ groups are biautomatic
(\cite{gersten}).  The interest here lies in the alternate method of
proof.  The following result is new.

\begin{thm}\label{generalization}

Let $\Gamma$ be a graph with every edge labeled by a positive integer
greater than $2$.  If there is a way of orienting the edges in the
graph so that neither of the graphs in Figure~\ref{bad graphs} appear
as subgraphs, then $K_\Gamma$ has a metric of non-positive curvature
and the Artin group $G_\Gamma$ is biautomatic.

\end{thm}

\begin{proof}

Since by Lemma~\ref{graphs} there are no embedded loops of length less
than $6$ in the link $L_\Gamma$, the metric which assigns a length of
$1$ to each $1$-cell and gives each $2$-cell the metric of a Euclidean
equilateral triangle will be non-positively-curved (\cite{ballmann}).
And since $K_\Gamma$ is an $A_2$ complex, by \cite{gersten} the
corresponding group is biautomatic.
\end{proof}

\begin{exmp}

Let $\Gamma$ be a graph in which all of the vertices have even degree
and assign large integers (each at least $3$) to each of its edges.
If $\Gamma$ can be embedded in the plane so that all of its short
embedded loops (length $3$ or $4$) actually bound components of the
complement of $\Gamma$, then the corresponding Artin group $G_\Gamma$
will be biautomatic.  The proof is an easy application of
Theorem~\ref{generalization}.  In particular, we can use the planarity
of the graph to orient the edges.  Assign an oriention to each of the
regions of the complement so that regions separated by an edge are
given opposite orientations.  The condition on the degrees of the
vertices guarantees that this is possible.  Now use these orientations
to induce the orientations of the edges.  The condition on the
embedded loops now guarrantees that neither of the forbidden subgraphs
can occur, and by Theorem~\ref{generalization}, the group $G_\Gamma$
is biautomatic.

\end{exmp}

\begin{rem}

Theorem~\ref{generalization} can be extended to include commutation
relations under the following restrictions.  Since the link of a
commutation relation includes special edges in {\em both} directions,
edges labeled by a $2$ in $\Gamma$ must be considered as oriented in
both directions.  In particular, when searching for a subgraph of type
$A$ or $B$, the edges labeled $2$ must be considered as `wildcards'
which can adopt either orientation as needed.  If there is an
orientation of the other edges of $\Gamma$ so that these subgraphs do
not appear, then the corresponding complex $K_\Gamma$ has a metric of
non-positive curvature and the group $G_\Gamma$ is biautomatic.

\end{rem}

Using various techniques, researchers have shown that other classes of
Artin groups are biautomatic. This is known, for example, for the
Artin groups of finite type (R.~Charney \cite{charney}) and for the
Artin groups of extra-large type (D.~Peifer \cite{peifer}).  Despite
this progress, the following conjecture remains open.

\begin{conj}

All Artin groups are biautomatic.

\end{conj}

\end{document}